\documentclass[a4paper,10pt]{article}
\usepackage{amssymb}
\usepackage[babel=true,kerning=true]{microtype}

\usepackage{amsthm,leqno}

\usepackage{amsfonts,amsmath,amssymb,enumerate}
\usepackage{graphicx}
\usepackage{tikz}
\usetikzlibrary{decorations.markings}
\usetikzlibrary{plotmarks}
\usetikzlibrary{patterns}

\newcommand{\ds}{\displaystyle}

\newcommand{\be}{\begin{equation*}}
\newcommand{\ee}{\end{equation*}}
\newcommand{\beq}{\begin{equation}}
\newcommand{\eeq}{\end{equation}}
\newcommand{\begincal}{\begin{eqnarray*}}
\newcommand{\fincal}{\end{eqnarray*}}

\newtheorem{thm}{Theorem}[section]

\newtheorem{prop}{Proposition}[section]
\newtheorem{defi}{Definition}[section]

\newcommand{\eps}{{\varepsilon}}
\newcommand{\R}{{\mathbb R}}
\newcommand{\D}{{\mathbb D}}
\newcommand{\C}{{\mathbb C}}
\newcommand{\N}{{\mathbb N}}
\newcommand{\h}{{\mathbb H}}
\newcommand{\Z}{{\mathbb Z}}
\newcommand{\Hy}{{\mathbb H}}

\def\eps{\varepsilon}
\def\ds{\displaystyle}

\title{Optimal estimate for the gradient of Green's function on degenerating surfaces and applications}
\author{Paul Laurain \& Tristan Rivi\`ere}

\begin{document}
\maketitle
\begin{abstract}
In this paper we prove a uniform  estimate for the gradient of the Green function on a closed Riemann surface, independent of its conformal class, and we derive compactness  results for immersions with $L^2$-bounded second fundamental form and for riemannian surfaces of uniformly bounded gaussian curvature entropy.
\end{abstract}
\noindent{\bf Math. Class.  32G15, 30F10, 53A05, 53A30, 35J35}

\section*{Introduction}
Let $\Sigma$ a closed smooth surface of genus $g$. We can endow $\Sigma$ with a metric $h$, then thanks to the uniformization theorem, see \cite{Jost06} or \cite{HPS}, there exists in the conformal class of $h$, i.e. the set of metric on $\Sigma$ which can be written $e^{2u} h$ where $u$ is a smooth function, a metric of constant curvature, equal to $1$ if $g=0$, $0$ if $g=1$ and  $-1$ otherwise. The sphere case is very particular, since the conformal group is not compact but this case is not of great interest here since there is only one conformal class. In the hyperbolic case the metric is unique and  in the torus it is also true up to normalized the area. In the following, we assume this normalization and we still denote by $h$ the metric of constant curvature (when $g\geq 1$) and we associate  to $h$ its Laplace-Beltrami operator $\Delta_{h} $. Then there exists, up to normalization, a unique nonnegative Green function $G_h$ associated to   $\Delta_{h} $.\\
  The main goal of this paper is to give estimates on $G_h$ independently of the conformal class defined by the metric $h$. This is a very classical subject in the theory of Riemann surfaces strongly related with the behavior of the spectrum of the Laplace operator, see \cite{Buser}.
Let $(\Sigma,h_k)$ a sequence of hyperbolic surface whose conformal class degenerate, that is to say that some geodesics are pinching. Let us assume that there is only one degenerating geodesic $\gamma_k$, let denote $\Sigma_\infty$ its nodal limit, see section 1 for precise definition, then Ji, see \cite{Ji},  proved that  If $\gamma_k$ does not separate $\Sigma_k$ then
$G_k$ is uniformly bounded  on every compact of $\Sigma_\infty \times \Sigma_\infty$, else
$\displaystyle \lim_{k\rightarrow +\infty} \vert G_k\vert = +\infty$ on  $\Sigma_\infty\times \Sigma_\infty$.\\

Here we see that we have a very different behavior with respect to the manner the conformal class degenerates. Is a similar behavior is possible for derivatives? Indeed, formally we can write
$$G_k(x,y) =\sum_{i\geq 1} \frac{\varphi_i^k(x)\varphi_i^k(y)}{\lambda_i^k},$$ \\ 
 where $\lambda_i^k$ and $\varphi_i^k$ are respectively the $i$th  (non vanishing) eigenvalue and the $i$th  (non constant) eigenfunction of $\Delta_{h_k}$, repeating indices according to multiplicity. Of course if the the nodal limit is disconnected then the first eigenvalue goes to zero while the first eigenfunction goes to a  positive locally constant function depending on the genus of each connected component. Looking at derivatives instead, one can expect a better behavior of the Green function, even in the collar region.  However, the gradient of the Green function gets a simple pole on the diagonal, hence it is not in $L^2$. The main result of this paper establishes that it is however true in a slightly weaker norm.
 
 \begin{thm}
\label{t1}
Let $\Sigma$ be a closed surface then there exits $C$ a positive constant and an integer $N$ depending only on the genus of $\Sigma$ such that for any metric $h$ on $\Sigma$ with constant curvature equal to $1$, $0$ or $-1$ and with normalized volume\footnote{$4\pi$ in the sphere case and $1$ in the torus case}, and any Green function associated to $h$, there exists a finite atlas of $N$ conformal charts $(U_i,\psi_i)$, such that for any $y\in\Sigma$ we get 
\beq
\label{mq}
 \sup_{t>0} t^2 \left\vert \left\{ x\in V_i \; \vert \; \vert d_x G_h^i (x, y ) \vert \geq t \right\}\right\vert \leq C,
 \eeq
 where $V_i= \psi_i(U_i)$ and $G_h^i(\, . \, ,y)=(\psi_i)_{*} (G_h(\, . \, , y))$.
\end{thm}

We can remark that on a fix Riemann surface $(\Sigma,h)$, the Green function is always bounded for this weak $L^2$-norm. We can deduce it from the standard pointwise estimate, see \cite{Aubin}, 
$$\vert d_x(G_h^i (x, y ))\vert\leq \frac{C_h}{d_h(x,y)} ,$$
but of course this estimate depends on the metric we take on $\Sigma$.\\

At the knowledge of authors, this result is the first control of the Green function independent of the conformal class. Moreover it looks quite optimal since the result is clearly false in $L^2$. In fact the atlas is very explicit, since, for intense considering the hyperbolic case, the surface divides in thick and thin part, on the thick part we can consider any disc with radius smaller than the injectivity radius. And in the thin part, using collar lemma, chart are given by degenerating annuli.\\

This result is optimal in the sense that we also prove that the weak $L^2$-norm of the Green function computed with respect to the intrinsic metric is not bounded when the singularity hold to a collapsing region, which is make clear  in the torus case by the  proposition \ref{t2}.\\

Regarding the proof of the theorem, once we have rule out the trivial case of the sphere then we treat the case of a degenerating torus and of an hyperbolic surface quite differently. For the torus, our proof relies on an estimate of the coefficient of the Fourier decomposition on a long thin cylinder using the periodicity condition, and in the hyperbolic case, it relies on the coarea formula and the decomposition of the surface in thin and thick part.\\

In the last sections, we give some applications of theorem \ref{t1} in differential geometry both from extrinsic and intrinsic point of view. First, we prove that the gradient of the conformal factor of an immersion with $L^2$-bounded second fundamental form is uniformly bounded in $L^{2,\infty}$, up to chose a convent atlas given by theorem \ref{t1}. Then we deduce a theorem of weak compactness for sequences of immersions with  $L^2$-bounded second fundamental. This last result was partially already proved by Kuwert and Li and the second author, see \cite{KL} and \cite{Riviere10}. Finally, in the last section we prove that considering a sequence of Riemann surface with bounded total curvature and entropy(see below for precise definition) then we can  find a finite conformal atlas in which the conformal factor is uniformly bounded.\\
  
\medskip {\bf Acknowledgements} : The first author was visiting  the {\it Forschungsinstituts f\"ur Mathematik} at E.T.H. (Zurich)  when this work started, he would like to thank it for its hospitality and the excellent working conditions.

\section{Preliminaries}

\subsection{Lorentz spaces}
Here we recall some classical facts about Lorentz spaces, \cite{Gra1} for details.\\

\begin{defi}  Let $D$ be a domain of $\R^k$, $p \in (1,+\infty)$ and  $q \in [1,+\infty]$. The Lorentz space $L^{p,q}(D)$ is the set of measurable functions
$f : D\rightarrow \R$ such that
$$\vert f\vert_{p,q}= \left(\int_{0}^{+\infty} \left( t^\frac{1}{p} f^{*}(t)\right)^q \frac{dt}{t}\right)^\frac{1}{q} <+\infty \hbox{ if } q<+\infty$$
or
$$\vert f\vert_{p,\infty}= \sup \left(t^\frac{1}{p} f^{*}(t)\right ) \hbox{ if } q=+\infty$$
where   $f^*$ the decreasing rearrangement of $f$.
\end{defi}

$\vert\, \, \vert_{p,q}$ happens to be a quasi norm  equivalent to a norm for which $L^{p,q}$ is a Banach space. Each $L^{p,q}$ may be seen as a deformation of $L^p$. For instance, we have
the strict inclusions
$$L^{p,1} \subset L^{p,q'} \subset L^{p,q''}\subset L^{p,\infty},$$
if $1 < q'< q''$. Moreover,
$$L^{p,p} = L^p.$$
Furthermore, if $\vert D \vert$ is finite, we have that for all $q$ and $q'$,
$$p > p' \Rightarrow  L^{p,q} \subset L^{p',q'}.$$

Using the fact that $f^*(t)=\inf\{s>0 \hbox{ s.t. } d_f(s)\leq t\}$ where $d_f$ is the distribution function of $f$, we see that the $L^{2,\infty}$ norm of $f$  is finite if and only if  $\displaystyle \sup_{t>0} t^2 \left\vert \left\{ x\in D \; \vert \; \vert f (x, \, . \, ) \vert \geq t \right\}\right\vert$ is finite.\\

Finally, for $p \in (1,+\infty)$ and  $q \in [1,+\infty]$, $L^{\frac{p}{p-1},\frac{q}{q-1}}$  is the dual of $L^{p,q}$.
\subsection{Degenerating Riemann surfaces}
Here we remind the Deligne-Mumford's description of the loss of compactness of the conformal class for a sequence of Riemann surfaces with fixed topology, see \cite{Hummel} for details.\\

Let $(\Sigma, c_k)$ a sequence of closed Riemann surface of fixed genus $g$. If $g=0$ then the conformal class is fixed since there is only one conformal class on the sphere. If $g=1$ then, we know that, $(\Sigma, c_k)$ is conformally equivalent to 
$\ds \R^2 / \left(\frac{1}{\sqrt{\Im(v_k)}}\Z \times \frac{v_l}{\sqrt{\Im( v_k)}}\Z\right)$ where $v_k$ lies in the fundamental domain $\{z\in \C \hbox{ s.t. } \vert\Re(z)\vert\leq 1\hbox{ and } \vert z \vert \geq 1\}$ of $\Hy/\mathrm{PSL}_2(\Z)$, and we say that $c_k$ degenerates if $\vert v_k\vert \rightarrow +\infty$. If $g\geq 1$, let $h_k$ the hyperbolic metric associated with $c_k$, then $(\Sigma, c_k)$ degenerates if there exits a closed geodesic whose length goes to zero. In that case, up to a subsequence, there exists
\begin{enumerate}
\item an integer $N\in \{ 1, \dots, 3g-3\}$,
\item a sequence $\mathcal{L}_k=\{ \Gamma_k^i\, ;\, i=1\dots N\}$ of finitely many pairwise disjoint simple closed geodesics of $(\Sigma,h_k)$ with length converging to zero,
\item a closed Riemann surfaces $(\overline{\Sigma},\overline{c})$,
\item a complete hyperbolic surface $(\widetilde{\Sigma},\widetilde{h})$ with $2N$ cups $\{(q^i_1,q_2^i)\, ; \, i=1\dots N\}$ such that $\widetilde{\Sigma}$ has been obtain topologically after removing the geodesic of $\mathcal{L}_k$ to $\Sigma$ and after closing each component of the boundary of $\Sigma\setminus \mathcal{L}_k$ by adding a puncture $q_l^i$ at each of these component. Moreover $\overline{\Sigma}$ is topologically equal to $\widetilde{\Sigma}$ and the complex structure defined by $\widetilde{h}$ on $\widetilde{\Sigma}\setminus \{q_l^i\}$ extends uniquely to $\overline{c}$. We can also equipped $\overline{\Sigma}$ with a metric $\overline{h}$ with constant curvature, but not necessarily hyperbolic since the genus of $\overline{\Sigma}$ can be lower than the one of $\Sigma$.

\end{enumerate}
$(\widetilde{\Sigma},\widetilde{h})$ is called the nodal surface of the covering sequence and $(\overline{\Sigma},\overline{c})$ is its renormalization. These objects are related, in the sense that, there exists a diffeomorphism $\psi_k: \widetilde{\Sigma}\setminus \{q_l^i\}  \rightarrow \Sigma\setminus \mathcal{L}_k$ such that $\widetilde{h}_k = \psi_k^* h_k$ converge in $C^{\infty}_{loc}$ topology to $\widetilde{h}$.
\section{Proof of theorem \ref{t1}}

Before starting the proof, we present a "baby case" illustrating the difficulty for getting some $L^{2,\infty}$-estimate for functions whose laplacian is in $L^1$ on a long thin cylinder. On a fix domain, such an estimate is a classical result, see theorem 3.3.6 of \cite{Helein}. Let us now consider the cylinder $C_l= \frac{1}{\sqrt{2\pi l}} \left(S^1\times\left[-\frac{l}{2},\frac{l}{2}\right]\right)$ which is identified with $S^1\times\left[-\frac{l}{2},\frac{l}{2}\right]$ endowed with the conformal metric $g=\frac{1}{2\pi l}  (d\theta^2 +dt^2)$. We set $u_{l}(t,\theta) = \frac{t^2}{4\pi l}$ which solves  $\Delta_g u_l =1$ . Then $\vert d u_l\vert_g = \frac{t}{\sqrt{2\pi l}}$ and we easily check that $\Vert d u_l \Vert_{L^{2,\infty}_g}\sim l$. While, considering the conformal  chart   $\psi_l : A_l \rightarrow C_l$ with
$A_l=\D \setminus B\left(0, e^{-l}\right)$  and $$\psi_l(\theta,r)= \left(\cos(\theta), \sin(\theta), \ln(r)+\frac{l}{2}\right).$$
Then, ${\overline{u}}_l=u_l\circ \psi_l = \frac{\left(\ln(r)+\frac{l}{2}\right)^2}{4\pi l}$ is uniformly bounded in $L^{2,\infty}$ {\bf with respect to the euclidean metric}. Indeed,
$$\vert \nabla  {\overline{u}}_l\vert = \left\vert \frac{\ln(r)+\frac{l}{2}}{2\pi l r} \right\vert \leq \frac{1}{r}.$$

This fact also illustrate that despite its closeness to the $L^2$-norm, the $L^{2,\infty}$-one is not conformally invariant\footnote{Although, it is invariant by dilation.}. This is one of the reason why we need to construct a specific conformal atlas.\\

{\bf All along the proof, for any given chart,unless otherwise stated, all the norms are computed with respect to the euclidean metric.}\\

\noindent{\it Proof of theorem \ref{t1} :}\\

{\bf The sphere case :}\\
Any sphere with constant curvature is conformal to the standard one, then  the chart are given by south and north stereographic projection composed by the conformal diffeomorphism. The Green function (up to a constant) is the one of the standard sphere, its gradient is clearly bounded in $L^{2,\infty}$, hence there is nothing to prove for theorem \ref{t1}.\\

{\bf The torus Case :}\\

Let $(\Sigma_l, g_{l})$ be a sequence of flat tori of volume $1$. Thanks to uniformization theorem, see \cite{Jost06}, we know that,  $(\Sigma_l, g_{l})$ is isometric to $\ds \R^2 / \left(\frac{1}{\sqrt{\Im(v_l)}}\Z \times \frac{v_l}{\sqrt{\Im( v_l)}}\Z\right)$ where $v_l$ lies in a fundamental domain of $\Hy/\mathrm{PSL}_2(\Z)$. Of course in the following we assume that the sequence degenerate, i.e. $\vert v_l\vert \rightarrow +\infty$, else the metric strongly converge and also the Green function.\\

 We are going to treat first the rectangular torus  (the case of $v_l\in i\R$) and we will explain how to deduce from it the general case. Up to some normalizations, our torus is isometric to  long and thin cylinder : $C_l=\frac{1}{\sqrt{2\pi l}}\left(S^1 \times \left[-\frac{l}{2},\frac{l}{2}\right]\right)$ with the standard  identification of its boundary components. Then this cylinder is conformal to the annular $A_l=\D \setminus B\left(0, e^{-l}\right)$  through the following diffeomorphism
$$\psi_l(\theta,r)= \left(\frac{\cos(\theta)}{\sqrt{2\pi l}}, \frac{sin(\theta)}{\sqrt{2\pi l}}, \frac{\ln(r)+\frac{l}{2}}{\sqrt{2\pi l}}\right).$$
Let $G_l$ be  the pull back of a Green function on $A_l$. It satisfies\footnote{ This equation must be understood in a weak sense and be tested against smooth function of $A_l$ whose composition with $\psi_l^{-1}$ extends to a smooth function on $\Sigma_l$.}
$$\Delta_z G_l(\, . \, , w)= \delta_w - \frac{1}{2\pi l r^2 } \hbox{ on } A_l ,$$
and 
$$G_l\left((\theta,e^{-l}),w \right)= G_l \left((\theta,1),w \right) \hbox{ and } e^{-l} \partial_r G_l\left((\theta,e^{-l}),w \right)= \partial_r G_l\left((\theta, 1),w \right) \hbox{ for all } \theta.$$

Then  we split $G_l$ in three parts: a singular part $s_l$, a diffusion part $u_l$ and an harmonic part $g_l$, as follows $G_l= s_l+u_l+g_l$ where 
$$s_l\left((\theta,r),w\right) =\left\{\begin{array}{c}\frac{1}{2\pi} \ln\left\vert r e^{i\theta} -w\right\vert+ \frac{1}{2\pi} \ln\left\vert r e^{i\theta} -  e^l w\right\vert \hbox{ if } \vert w \vert \leq\frac{1}{2} \\  

\\
\frac{1}{2\pi} \ln\left\vert r e^{i\theta} -w\right\vert+ \frac{1}{2\pi} \ln\left\vert r e^{i\theta} -  e^{-l} w\right\vert \hbox{ if } \vert w \vert >\frac{1}{2}\end{array}\right. $$
and
$$u_l\left((\theta,r),w\right)=-\frac{1}{4\pi l} \left(\ln(r)\right)^2.$$

We easily check that, on the one hand  $\Delta s_l =\delta_w$ on $A_l$ and $\Vert \nabla s_l\Vert_{2,\infty}=O(1)$, independently of $l$ and $w$, and on  the other hand  $\Delta u_l =-\frac{1}{2\pi l r^2 }$ on $A_l$ and $\Vert \nabla u_l\Vert_{2,\infty}= O(1)$, in fact we get even the more precise estimate $ \vert \nabla u_l\vert= O\left(\frac{1}{r}\right)$. Finally we estimate $g_l$, in that aim we assume that $w\leq \frac{1}{2}$, the other case can be done in a similar way. Then $g_l$ satisfies
$$\Delta g_l= 0,$$
\beq
\label{e1}
g_l\left((\theta,1),w \right)-g_l\left((\theta,e^{-l}),w \right) = -\frac{3l}{4\pi}+\frac{1}{2\pi} \ln \left\vert \frac{e^{-l}e^{i\theta}-e^{l}w }{e^{i\theta}-w}\right\vert=  -\frac{3l}{4\pi}+ F_l(\theta),
\eeq
and  
\beq
\label{e2}
\begin{split}
\partial_rg_l\left((\theta,1),w \right)- e^{-l}\partial_r g_l\left((\theta,e^{-l}),w \right)& = \frac{1}{2\pi} \left(- \frac{\langle e^{i\theta}, e^{i\theta}-w\rangle }{\vert e^{i\theta}-w\vert^2} +  \frac{\langle e^{-l}e^{i\theta}, e^{-l}e^{i\theta}- e^lw\rangle }{\vert e^{-l}e^{i\theta}-e^l w\vert^2} + 1\right)\\
&=H_l(\theta).
\end{split}
\eeq

Using Fourier analysis, we can decompose $g_l$  as follows
$$g_l((\theta,r), w)= c_0+a_0 \ln(r)+\frac{1}{\sqrt{2\pi}}\sum_{n\in \Z} (a_n r^n+b_n r^{-n})e^{in\theta} .$$

Thanks to (\ref{e1}) we easily check that,
$$a_0 = O(1) .$$
On the one hand, thanks to (\ref{e1}) and (\ref{e2}), we get

\beq
\label{e3}
a_n (1-e^{-nl}) +b_n  (1-e^{nl})=\frac{1}{\sqrt{2\pi}} \int_0^{2\pi} F_l(\theta)e^{-in\theta}\, d\theta
\eeq
and 
\beq
\label{e4}
a_n  (1-e^{-nl}) -b_n  (1-e^{nl})=\frac{1}{n\sqrt{2\pi}} \int_0^{2\pi} H_l(\theta)e^{-in\theta}\, d\theta .
\eeq
On the other hand, 
\be
\begin{split}
\left\Vert \nabla \left(g_l -a_0 \ln(r) \right)\right\Vert_2^2 &=O\left( \sum_{n\in \Z} n^2 a_n^2 \int_{e^{-l}}^{1} r^{2n-1}\, dr + n^2 b_n^2 \int_{e^{-l}}^{1} r^{-2n-1}\, dt \right)\\
&=O\left(  \sum_{n\in \Z} n a_n^2 \left(1-e^{-2nl} \right)+ n b_n^2 \left(1-e^{2nl} \right) \right)
\end{split}
\ee
But thanks to  (\ref{e3}) and  (\ref{e4}) and the fact that $F_l$ and $H_l$ converge in $C^2(S^1)$, as $l$ goes to infinity, we have

$$ \vert a_n\vert = O\left(\frac{1}{n^2(1-e^{-nl})}\right) \hbox{ and }  \vert b_n\vert = O\left(\frac{1}{n^2 (1-e^{nl})}\right) $$
uniformly with respect to  $l$. Hence  

$$\sum_{n\in \Z} n a_n^2 \left(1-e^{-2nl} \right)+ n b_n^2 \left(1-e^{2nl} \right) =O(1)$$
uniformly with respect to $l$. Which prove that 

$$ \left\Vert \nabla \left(g_l - a_0 \ln(r) \right)\right\Vert_2^2  =O(1)$$

Finally we conclude that 
$$\Vert \nabla G_l \Vert_{2,\infty} = O(1),$$

This achieves the proof of theorem \ref{t1} in the case $v_l \in i \R$. In the general case, the torus is isometric to $(A_{l_k},g_k)$ where $l_k=2\pi \Im(v_k)$ and $g_k =\frac{1}{2\pi l_k r^2}\left( r^2d\theta^2 + \frac{dr^2}{\cos(\alpha_k)}\right)$ with $\alpha_k =\frac{\pi}{2} -\arg(v_k)\rightarrow 0$. Hence $G_{l_k}$ split as follows $G_{l_k} = s_k+u_k+g_k$ where 

$$s_k\left((\theta,r),w\right) =\left\{\begin{array}{c}\frac{1}{2\pi \cos(\alpha_k)} \ln\left\vert r e^{i\theta} -w\right\vert_{g_k} + \frac{1}{2\pi \cos(\alpha_k)} \ln\left\vert r e^{i\theta} -  e^l w\right\vert_{h_k} \hbox{ if } \vert w \vert \leq\frac{1}{2} \\  

\\
\frac{1}{2\pi \cos(\alpha_k)} \ln\left\vert r e^{i\theta} -w\right\vert_{g_k}+ \frac{1}{2\pi \cos(\alpha_k)} \ln\left\vert r e^{i\theta} -  e^{-l} w\right\vert_{g_k} \hbox{ if } \vert w \vert >\frac{1}{2}\end{array}\right. $$\\

and 

$$u_k\left((\theta,r),w\right)=-\frac{1}{4\pi l_k \cos(\alpha_k)} \left(\ln(r)\right)^2,$$
and then the analysis of $g_k$ is the same.\\

The fact that we cannot bound the weak-$L^2$ norm with respect to the intrinsic metric of the torus is illustrated by the following proposition.

\begin{prop} 
\label{t2}
There exists a sequence of metric $h_k$ on $\mathbb{T}^2$ with constant curvature and volume equal to $1$ which is unbounded in the moduli space and such that, for any $y\in\Sigma$, we get 

\beq
 \sup_{t>0} t^2 \left\vert \left\{ x\in \Sigma \; \vert \; \vert dG_{h_k} (x, y ) \vert \geq t \right\}\right\vert \rightarrow + \infty.
 \eeq
\end{prop}

\noindent{\it Proof of proposition \ref{t2}:}\\

In order to prove this proposition ,we go back to the case of degenerating rectangular tori viewed as a long thin cylinder: $C_l=\frac{1}{\sqrt{2\pi l}}\left(S^1 \times \left[-\frac{l}{2},\frac{l}{2}\right]\right)$ with the two ends being identified in an obvious way. Let $g_\phi(\theta,t) =G((\theta,t),(\phi,0))$, then if this function is uniformly bounded in $L^{2,\infty}$  it would be the case for $g(\theta,t) =\frac{1}{2\pi} \int_0^{2\pi} g_\phi(\theta,t) \; d\phi$, thanks to the invariant by rotation.   
Then we easily check that $g(\theta,t)=\frac{\vert t \vert}{4\pi}$, hence $\vert dg \vert = \frac{\sqrt{2\pi l}}{4\pi}$, here the norm is computed with respect to the metric $h =\frac{1}{2\pi l} (d\theta^2 + dt^2)$. Finally we easily check that  $\Vert d g \Vert_{L^{2,\infty}_h}\sim l$, which is a contradiction and prove proposition \ref{t2}. \hfill$\square$\\

{\bf The case of genus $\geq2$ :}\\
 
Let $(\Sigma, c_k)$ be  a sequence of Riemann surfaces of fixed  genus $g\geq 2$. Thanks to the uniformization theorem, see \cite{Jost06}, we know that, we can endow $\Sigma$ with a conformal metric $h_k$ such that   $(\Sigma, h_{k})$ is isometric to $\ds \h / \Gamma_k$ where $\Gamma_k$ is a discrete group of  $\mathrm{PSL}_2(\R)$.
Then, the Green functions associated to $h_k$ satisfy 
$$ (\Delta_{h_k})_x G_k(\, .\, , y) = \delta_y -\frac{1}{v_k} ,$$
where $v_k$ is the volume of $(\Sigma,h_k)$ which depends only on the genus thanks to Gauss-Bonnet theorem. In order to study the behavior of the Green's function and following the classic description of hyperbolic surfaces, see \cite{Hummel} or \cite{Hubbard}, we set $\delta<arcsinh(1)$ and  then we  split $(\Sigma, h_{k})$ in two parts: a thick part $E_k^\delta=\{ s\in \Sigma \, \vert \, \mathrm{injrad}((\Sigma,h_k),s)\geq\delta\}$ and a thin part $F_k^\delta =\{ s\in \Sigma \, \vert \, \mathrm{injrad}((\Sigma,h_k),s)< \delta\}$. Thanks to the decomposition theorem of Deligne-Mumford, we know that the sequence of metrics converges strongly on the thick part  and develops collar in the thin part. We are going to split our proof depending on the case whether  $y$ lies in the thick or int thin part. But before we prove a general estimate for Green functions on a closed surfaces.\\

{\bf Step 1: Coarea formula for Green's functions.}\\

In this step, $G_k$ is any Green function associated to $h_k$, there is no normalization. Volumes and lengths are taken with respect to $h_k$ here.\\

Let $t>0$, integrating by part on a level set, we get

$$\int_{G_k(x,y)=t} \frac{\partial G_k(x,y)}{\partial \nu_k}\, d\sigma_k(x)=\int_{G_k(x,y)\geq t} (\Delta_{h_k})_x G_k(x,y)\, dv_k(x)=1-\frac{vol(\{G_k(x,y)\geq t\} )}{v_k},$$
where $\nu_k$ is the exterior normal of the open set $\{ x\in \Sigma \; \vert \; G_k(x,y)> t \} $.
Then, we get
$$\int_{G_k(x,y)=t} |dG_k(x,y) |^2_{h_k}\, d\sigma_k(x) \leq 2 .$$
Let $1<p<2$ and $a>0$, then thanks to coarea-formula, see \cite{Simon} or \cite{Chavel06}, we get that
$$\int_{G_k(x,y)\geq a} |dG_k^{1-p/2}(x,y) |^2_{h_k} \, dv_k(x)=\int_a^{+\infty}\left(\int_{G_k(x,y)=t} \frac{|dG_k(x,y) |_{h_k}}{G_k^p(x,y)} \, d\sigma_k(x)\right)\, dt\leq  2\int_a^{+\infty}\frac{1}{t^p}\, dt\leq C_{p,a},$$
where $C_{p,a}$ is a positive constant depending only on $p$ and $a$. Moreover, we can also prove, considering negative level set, that
$$\int_{G_k(x,y)\leq -a} \vert d G_k^{1-\frac{p}{2}}(x,y)\vert^2_{h_k} \, dv_k(x)\leq C_{p,a},$$
where $C_{p,a}$ is a positive constant depending only on $p$ and $a$. Finally we get, 
\beq
\label{coarea}
\int_{\vert G_k(x,y) \vert \geq a} \vert d G_k^{1-\frac{p}{2}}(x,y)\vert^2_{h_k} \, dv_k(x)\leq C_{p,a},
\eeq
where $C_{p,a}$ is a positive constant depending only on $p$ and $a$.\\

{\bf Step 2: Estimate in the thick part.}\\

In order to obtain the estimate on the whole thick part, we will cover it by a finite number of balls with radius $\frac{\delta}{2}$, where have been chosen such that $0 < \delta< acrsinh(1)$. Since we consider a general sequence of Green functions $G_k(\, .\, , y_k)$, we have to pay attention to the location of these balls with respect to the singularity  $y_k$. In fact if $y_k$  is in the thick part then we will center one of the ball of the covering  at $y_k$ and then the others won't have to deal with this singularity.\\

Let $x_k\in E_k^\delta$ and we first assume that $y_k \not\in B_{h_k}\left(x_k,\frac{\delta}{2}\right)$. Then $B_{h_k}\left(x_k,\frac{\delta}{2}\right)$ is isometric to $B\left(0, \tanh\left( \frac{\delta}{4}\right)\right)$ in the Poincar\'e disc. In the following, we make all computations in the conformal chart $B(0,3r)$ with $r=\frac{\tanh\left( \frac{\delta}{4}\right)}{3}$ and the metric $h_p=\frac{4dx^2}{(1-\vert x\vert^2)^2}$. But the hyperbolic metric is equivalent to the euclidean one on this ball.\\

On $B(0,3r)$ we decompose $G_k$ as follows
$$ G_k(\, .\, , y_k)= u_k + g_k,$$

$$ G_k(\, .\, , y_k)= u_k + g_k,$$
where $u_k(x)= \frac{1}{v_k} \ln\left(\frac{1}{ 1-\vert z\vert^2}\right)$ and $g_k$ be a smooth harmonic function. Hence we can apply (\ref{coarea}), with $p=\frac{3}{2}$ and $a=\frac{2}{v_k} \ln\left(\frac{1}{ 1-(3r)^2}\right)$ to $\widetilde{G}_k=G_k - g_k(0)$, which gives that 
$$\int_{\{\vert \widetilde{G}_k(x,y_k)\vert \geq a\}\cap B(0,3r) } \vert d \widetilde{G}_k^{\frac{1}{4}}(x,y_k)\vert^2_{h_p} \, dv_{h_p}(x)\leq C,$$
where $C$ is a positive constant depending only on the genus and $\delta$. Then, by the mean value property, there exists $\rho\in[2r,3r]$ such that
$$\int_{\{\vert \widetilde{G}_k(x_k,y) \vert \geq a\}\cap \partial B(0,\rho)} \vert d \widetilde{G}_k^{\frac{1}{4}}(x,y_k)\vert_{h_p} \, d\sigma_{h_p}(x) \leq C,$$
where $C$ is positive constant depending only on the genus and $\delta$. Then, using the fact harmonic functions satisfy the mean value property, we get  that the mean value of $g_k-g_k(0)$ is $0$ on $ \partial B(0,\rho)$ and we easily deduce that
$$ \vert G_k(x,y_k)-g_k(0)\vert\leq C, \hbox{ for all } x \in  \partial B(0,\rho),$$
where $C$ is positive constant depending only on the genus and  $\delta$. Then using classical elliptic estimate, we have
\beq
\label{ei}
\Vert d G_k(\, .\, , y_k)\Vert_{L^{\infty}_{h_p}(B(0,r))} \leq C, 
\eeq
where $C$ is positive constant depending only on the genus and $\delta$.  Then, in the case $x_k=y_k$, equivalently the ball we consider is centered at the singularity of $G (x,\,. \, )$, we obtain the same estimate decomposing $G_k$ as follows
$$ G_k(\, .\,, y_k)= s_k+u_k + g_k,$$
where $s_k(x)=\frac{1}{2\pi}\ln(\vert x\vert)$. But, of course, in that case, due to the presence of $s_k$, the estimate is in $L^{2,\infty}$. Finally covering the thick part with a uniformly bounded number of balls we get the desired estimate on $G(\, . \, ,y_k)$ on the thick part. Indeed, either $y_k$ is not in the thick part and the result will follows directly from (\ref{ei}), or we start by taking a ball centered at $y_k$ and then we cover the rest of the thick part by balls which does not contain $y_k$.\\

{\bf Step  3: Estimate in the thin part.}\\

Let $x_k\in F_k^\delta$ and $y_k \in \Sigma$ two converging sequences in $\overline{\Sigma}$. First, thanks to the collar lemma, see \cite{Hubbard}, we know that each connected component of the thin part (i.e. at most 3g-3), contains a simple closed geodesic $\gamma_k$ of length $\eps_k=l(\gamma_k)< 2  arcsinh(1)$, and is isometric to 

$$B_k=\left\{ z=re^{i\varphi}\in \h : 1\leq r \leq e^{\eps_k}, arctan\left(sinh\left(\frac{\eps_k}{2}\right)\right)<\varphi<\pi - arctan\left(sinh\left(\frac{\eps_k}{2}\right)\right)\right\},$$
where the geodesic corresponds to $\left\{ r e^{i\frac{\pi}{2}}\in \h : 1 \leq r \leq e^{\eps_k}\right\}$ and the line $\{r=1\}$ and $\{ r=e^{\eps_k}\}$ are identified via $z\mapsto e^{\eps_k} z$. It is often easier to consider the following cylindrical parametrization. Let $\varphi_k= arctan\left(sinh\left(\frac{\eps_k}{2}\right)\right)$ and we set
$$C_k=\left\{ (\cos(\theta), sin(\theta),t) \, \vert \, 0\leq \theta <2\pi ,  \frac{2\pi}{\eps_k}\varphi_k <t<\frac{2\pi}{\eps_k} \left(\pi - \varphi_k\right)\right\}$$
equipped with the metric 
$$h_c=\left(\frac{\eps_k}{2\pi sin\left(\frac{\eps_k t}{2\pi}\right) }\right)^2 (d\theta^2+ dt^2 ),$$
where the geodesic correspond to $\left\{ t=\frac{\pi^2}{\eps_k}\right\}$.\\

We are going to make the proof assuming that $y_k$ lies in the thin part. When this is not the case the proof carries over after the simply operation consisting of withdrawing the singular part $s_k$. We can also assume that $y_k\not\in \left(\left[\frac{2\pi}{\eps_k}\varphi_k, \frac{2\pi}{\eps_k}\varphi_k +\frac{\delta}{10}\right]\cup\left[\frac{2\pi}{\eps_k} \left(\pi - \varphi_k\right)-\frac{\delta}{10},\frac{2\pi}{\eps_k} \left(\pi - \varphi_k\right)\right]\right) \times S^1$, replacing $\delta$ by $\frac{\delta}{2}$ if necessary.\\

Then, as for the torus case, we choose an annulus as conformal chart. Precisely, let $A_k= \D \setminus B(0,e^{-l_k})$ and  $\psi_k:A_k \rightarrow C_k$ defined as follows
$$\psi_k(\theta,r)= \left(\cos(\theta), sin(\theta), \ln(r)+\frac{2\pi}{\eps_k} \left(\pi - \varphi_k\right) \right),$$

where $l_k =\frac{2\pi}{\eps_k} \left(\pi -2 \varphi_k\right)$. Then, the pull back of a Green function on $A_k$, that we keep denoting $G_k$, satisfies

$$\Delta_z G_k(\, . \, , w_k)= \delta_{w_k} - \left(\frac{\eps_k}{r2\pi \sin\left(\frac{\eps_k}{2\pi} \left(\ln(r)+\frac{2\pi}{\eps_k} \left(\pi - \varphi_k\right) \right)\right) }\right)^2 \hbox{ on } A_k ,$$

where $\psi_k(w_k)=y_k$ with $w_k\in B(0,e^{-\frac{\delta}{10}})\setminus B(0,e^{-l_k+\frac{\delta}{10}})$ .\\

First of all, thanks to our previous step, see (\ref{ei}), we remark that
\beq
\label{ei21}
\vert \nabla G_k(\, .\,, w_k) \vert \leq  \frac{C}{r}\hbox{ on } B(0, e^{-l_k + \frac{\delta}{10}})\setminus B(0, e^{-l_k}),
\eeq
and 

\beq
\label{ei22}
\vert \nabla G_k(\, .\,, w_k) \vert \leq C \hbox{ on } B(0, 1)\setminus B(0, e^{-\frac{\delta}{10}}),
\eeq

where $C$ is a positive constant depending only on the genus and $\delta$. Then we split $G_k$ as follows 
$$G_k((\theta,r),w_k)=s_k(\theta,r)+ u_k(\theta,r)+ g_k(\theta,r), $$
where 
$$u_k(\theta,r)= \frac{\ln\left( \sin\left(\frac{\eps_k }{2\pi} \left(\ln(r) +\frac{2\pi}{\eps_k}\right)\right)\right)}{v_k},$$
and 
$$s_k(\theta,r)= \frac{1}{2\pi}\ln\left(\left\vert re^{i\theta}-w_k\right\vert\right).$$
We easily check that 
$$ \Delta u_k=  - \left(\frac{\eps_k}{r2\pi \sin\left(\frac{\eps_k}{2\pi} \left(\ln(r)+\frac{2\pi}{\eps_k} \left(\pi - \varphi_k\right) \right)\right) }\right)^2 $$
\beq
\label{ll1}
\begin{split}
 \Vert \nabla u_k \Vert_2^2 &= \frac{1}{2\pi v_k^2} \int_{e^{-l_k}}^{1} \left(\frac{\cos\left(\frac{\eps_k}{2\pi} \left(\ln(r)+\frac{2\pi}{\eps_k} \left(\pi - \varphi_k\right) \right)\right)}{\sin\left(\frac{\eps_k}{2\pi} \left(\ln(r)+\frac{2\pi}{\eps_k} \left(\pi - \varphi_k\right) \right)\right)} \right)^2 \frac{\eps^2_k}{r} dr   \\
 & 
=  \frac{1}{2\pi v_k^2} \int_{e^{-l_k}}^{1} \left( \frac{-\eps^2_k}{r}  +\frac{1}{\left(\sin\left(\frac{\eps_k}{2\pi} \left(\ln(r)+\frac{2\pi}{\eps_k} \left(\pi - \varphi_k\right) \right)\right)\right)^2 } \frac{\eps^2_k}{r}\right) dr 
 \\
&
=\frac{1}{ v_k^2} \left[ -\frac{\eps^2_k}{2\pi} \ln(r) -\eps_k  \frac{\cos\left(\frac{\eps_k}{2\pi} \left(\ln(r)+\frac{2\pi}{\eps_k} \left(\pi - \varphi_k\right) \right)\right)}{\sin\left(\frac{\eps_k}{2\pi} \left(\ln(r)+\frac{2\pi}{\eps_k} \left(\pi - \varphi_k\right) \right)\right)} \right]_{e^-{l_k}}^1 \leq C,
\end{split}
\eeq

$$ \Delta s_k = \delta_{w_k},$$
and 
\beq 
\label{ll2}
\Vert \nabla s_k \Vert_{2,\infty} \leq C,
\eeq
where $C$ is a positive function depending only on the genus and $\delta$. 

Then $g_k$, which has been obtained from $G_k$ after subtracting $s_k$ and $u_k$, is a smooth harmonic function. Let $\overline{g}_k(r)$ be the mean value of $g_k$ on the circle of radius $r$ centered at $0$. It is also harmonic and radial, hence 
$\overline{g}_k(t)=a_k \ln(r)+ b_k$. Moreover, thanks to (\ref{ei21}) and (\ref{ei22}), we get that 

\beq
\label{ei31}
\vert \nabla g_k  \vert \leq  \frac{C}{r}\hbox{ on } B(0, e^{-l_k + \frac{\delta}{20}})\setminus B(0, e^{-l_k}),
\eeq
and 

\beq
\label{ei32}
\vert \nabla g_k \vert \leq C \hbox{ on } B(0, 1)\setminus B(0, e^{-\frac{\delta}{20}}).
\eeq

In particular, $a_k$ is uniformly bounded and we get,

\beq
\label{ei4}
\Vert \nabla \overline{g}_k \Vert_{L^{2,\infty}\left(B(0,1)\setminus B(0, e^{-l_k}) \right)} \leq C,
\eeq

Then using the fact the mean value of $g_k-\overline{g}_k$ is zero and the previous estimate, we get that
\beq
\Vert g_k-\overline{g}_k \Vert_{L^\infty\left(\left(B(0, e^{-l_k + \frac{\delta}{20}})\setminus B(0, e^{-l_k})\right) \cup \left( B(0, 1)\setminus B(0, e^{-\frac{\delta}{20}}) \right)\right)} \leq C,
\eeq
where $C$ is positive constant depending only on the genus and $\delta$. Then, since $g_k-\overline{g}_k$ is harmonic and with radial mean value equal to zero, 
\beq
\label{ei3}
\Vert \nabla  (g_k -\overline{g}_k) \Vert_{L^{2}\left(B(0, e^{-\frac{\delta}{10}}) \setminus B(0, e^{-l_k + \frac{\delta}{10}}) \right)} \leq C,
\eeq
where $C$ is positive constant depending only on the genus and $\delta$. The last inequality can be proved using the furrier decomposition and remarking that $g_k-\overline{g}_k$ has no logarithmic part. Finally, thanks to (\ref{ll1}), (\ref{ll2}), (\ref{ei4}) and (\ref{ei3}), we get the desired estimate, which concludes the proof of the theorem \ref{t1}.\hfill$\blacksquare$

\section{Weak compactness result for immersions with second fundamental form bounded in $L^2$ }

The first application of theorem \ref{t1} regards the control of the conformal factor for immersions with $L^2$-bounded second fundamental form. Before to state the main result, we shall first remind the notion of weak immersions introduced by the second author in \cite{Riviere10} .\\

Let $\Sigma$ a smooth compact surface equipped with a reference smooth metric $g_0$. One define the Sobolev spaces $W^{k,p}(\Sigma,\R^m)$ of measurable maps from $\Sigma$ into $\R^m$ into the following way
$$W^{k,p}(\Sigma,\R^m)=\left\{ f:\Sigma\rightarrow \R^m  \mathrm{ measurable s.t. } \sum_{l=0}^k \int_\Sigma \vert \nabla^l f\vert^p_{g_0} \, dv_{g_0} < +\infty\right\} .$$
Since $\Sigma$ is compact it is not difficult to see that this space is independent of the choice we have made of $g_0$.\\

Let $f\in W^{1,\infty}(\Sigma, \R^m)$, we define $g_f$ to be the following symmetric bilinear form
$$g_f(X,Y)= \langle df(X), df(Y)\rangle,$$
and we shall assume that there exists $C_f>1$ such that 
\beq
\label{C}
C_f^{-1} g_0(X,X) \leq g(X,X) \leq C_f g_0(X,X).
\eeq
For such a map, we can define the Gauss map as being the following measurable map in $L^\infty(\Sigma)$ taking values int the Grassmannian of oriented $m-2$-planes of $\R^m$,
$$\vec{n}_f = \star \frac{\frac{\partial f}{\partial x}\wedge \frac{\partial f}{\partial x}}{\left\vert \frac{\partial f}{\partial x}\wedge \frac{\partial f}{\partial x}\right\vert}.$$
We then introduce the space $\mathcal{E}_\Sigma$ of weak immersions of $\Sigma$ with bounded second fundamental form as follow:

$$ \mathcal{E}_\Sigma = \left\{\begin{array}{c}\Phi \in W^{1,\infty} (\Sigma) \hbox{ which satisfies (\ref{C}) for some }C_\Phi>0  \\
\\ \
\mathrm{ and } \int_\Sigma \vert  d\vec{n}_{\Phi}\vert_g^2 \, dv_{g} < +\infty \end{array}\right\},$$
where $g=\Phi^* \xi$.\\

It is proved in \cite{NCRiviere} that any weak immersion defines a smooth conformal structure on $\Sigma$. Let $\Phi\in \mathcal{E}_\Sigma$, we denote by $\pi_{\vec{n}_\Phi}$ the orthonormal projection of vector in $\R^m$ onto the $m-2$-plane given by $\vec{n}_\Phi$. With these notations the second fundamental form of the immersion at $p$ is given by
$$\forall X,Y \in T_p \Sigma \ \ \vec{\mathbb{I}}_p(X,Y)=\pi_{\vec{n}_\Phi} d^2 \Phi(X,Y),$$
and the mean curvature vector of the immersion at $p$ is given by 
$$ \vec{H}=\frac{1}{2} \mathrm{tr}_g(\vec{\mathbb{I}}).$$
A natural quantity while considering such immersions is the Lagrangian given by the $L^2$-norm of the second fundamental form :
$$E(\phi)=\int_\Sigma \vert \vec{\mathbb{I}}\vert_g^2 \, dv_g.$$
An elementary computation, using Gauss-Bonnet formula, gives
$$E(\phi)=\int_\Sigma \vert \vec{\mathbb{I}}\vert_g^2 \, dv_g=\int_\Sigma \vert  d\vec{n}_{\Phi}\vert_g^2 \, dv_{g} = 4 W(\phi)-4\pi \chi(\Sigma), $$
where $\chi(\Sigma)$ is the Euler characteristic and 
$$W(\Phi)= \int_\Sigma \vert \vec{H}\vert_g^2 \, dv_g,$$
is the so called Willmore energy. 

\begin{thm}
\label{a1}
 Let $(\Sigma,c_k)$ be a sequence of closed Riemann surface of fixed genus greater than one. Let denote $h_k$ the metric with constant curvature (and volume equal to one in the torus case) in $c_k$ and $\Phi_k$ a sequence of weak conformal immersion of $\Sigma$ into $\R^m$, i.e.
$$\Phi_k^*\xi =e^{2u_k}h_k,$$
where $u_k\in L^\infty(\Sigma)$. Then there exists a finite conformal atlas $(U_i,\psi_i)$ independent of $k$ and a positive constant $C$ depending only on the genus of $\Sigma$, such that 
$$\Vert d v^i_k \Vert_{L^{2,\infty}(V_i)}\leq C W(\Phi_k),$$
where $v^i_k$ is the conformal factor of $\Phi^k \circ \psi_i^{-1}$ in $V_i=\psi_i(U_i)$, i.e. $v_k^i=\frac{1}{2} \ln \left\vert \frac{\partial \Phi^k\circ \psi_i^{-1}}{\partial x} \right\vert=\frac{1}{2} \ln \left\vert \frac{\partial \Phi^k\circ \psi_i^{-1}}{\partial y}\right\vert$.
\end{thm}

\noindent{\it Proof of theorem \ref{a1}:}\\

Let $K_g= 0 \hbox{ or} -1$ if the genus $g$ of $\Sigma$ is $ 1$ or greater than $1$ and be $K_{\Phi_k^*\xi}$ the Gauss curvature associated to $\Phi_k^*\xi$. It is classical that $u_k$ satisfies the following Liouville equation
\beq
\label{eq}
-\Delta_{h_k} u_k= K_{\Phi_k^*\xi }  e^{2u_k}-K_g.
\eeq
Let $G_k$ be the nonnegative Green function of $(\Sigma,h_k)$, then using the representation formula, we get that
\beq
\label{conv}
 u_k = G_k \star \left(K_{\Phi_k^*\xi } e^{2u_k}-K_g\right) .
 \eeq
We have the following straightforward estimate,
\beq
\label{bl1}
\int_\Sigma \vert K_{\Phi_k^*\xi }  e^{2u_k}\vert\, dv_{h_k}=\int_\Sigma \vert K_{\Phi_k^*\xi }\vert \, dv_{\Phi_k^*\xi } \leq \frac{1}{2} \int_\Sigma \vert \vec{\mathbb{I}}_{\Phi_k^*\xi }\vert \, dv_{\Phi_k^*(\xi) } \leq W(\Phi_k) ,
\eeq
this proves that the right hand side of (\ref{eq})  is uniformly bounded in $L^1$-norm by $W(\Phi_k)$ with respect to the metric $h_k$. Then let $(U_i,\psi_i)$ the conformal atlas given by theorem \ref{t1} and let $\phi_k^i=\Phi_k \circ \psi_i^{-1} : V_i \rightarrow \R^m$ and $v_k^i :V_i \rightarrow \R$ such that 
$(\phi_k^i)^*(\xi)= e^{2 v_k^i} dz^2$. First we observe that $v_k^i = u_k \circ \psi_i^{-1} + w_k^i$ where 
$w_k^i :V_i \rightarrow \R$ such that  $(\psi_i)^*(h_k)= e^{2 w_k^i} dz^2$. Moreover, we can easily check that  $\nabla w_k^i$ is uniformly bounded in $L^{2,\infty}$ since in the torus case a chart is given by an annulus and $w_k^i= - \ln(r)+ c_k$ and in the hyperbolic case the chart is either a disc with radius strictly less than $1$ and  $w_k^i= - \ln(1-r^2)+ c_k^i$ or an annulus $A_k= \D \setminus B(0,e^{-l_k})$ where $l_k =\frac{2\pi}{\eps_k} \left(\pi -2 \varphi_k\right)$ with $w_k^i=  -\ln\left(r \sin\left(\frac{\eps_k}{2\pi} \left(\ln(r)+\frac{2\pi}{\eps_k} \left(\pi - \varphi_k\right) \right)\right)\right)+ c_k^i$.\\

Then it suffices to check that $\nabla (u_k \circ \psi_i^{-1})$ is uniformly bounded in $L^{2,\infty}(V_i)$ in order to prove the theorem. Thanks to (\ref{conv}), we have
\beq
u_k(y) =\int_\Sigma G_k(x,y) F_k(x)  \, dv_{h_k}  + \frac{1}{\mathrm{vol}( \Sigma ,h_k)}\int_\Sigma u_k \,dv_{h_k}  ,
\eeq 
where $F_k= K_{\Phi_k^*\xi }  e^{2u_k}-K_g$. Hence, 

\beq
\nabla_y u_k\circ \psi_i^{-1} (y) =\int_\Sigma \nabla_y G_k(x, \psi_i^{-1}(y)) F_k(x)  \, dv_{h_k}   ,
\eeq 

then using the fact $ \sup_{x\in \Sigma} \Vert \nabla_y G_k(x,\psi_i^{-1}(y)) \Vert_{L^{2,\infty}}$ is uniformly bounded, thanks to theorem \ref{t1}, thanks to (\ref{bl1}) and standard inequality on convolution we get that $\Vert \nabla (u_k\circ \psi_i^{-1} )\Vert_{L^{2,\infty}(V_i)}$ is uniformly bounded which concludes the proof of theorem \ref{a1}.\hfill$\blacksquare$\\

Then our second application concerns the weak compactness of the conformal immersion with $L^2$-bounded second fundamental form. The following result was proved first in \cite{Riviere10}
when the conformal classes of the surfaces do not degenerate and has been extended to the general case
of degenerating riemmann surfaces in \cite{KL}. We shall present a different approach for proving this result as being a consequence of our main theorem \ref{t1}

\begin{thm} 
\label{a2}
Let $\Sigma$ a  closed surface of genus strictly greater than $1$ and $\Phi_k\in \mathcal{E}_\Sigma$ a sequence of weak immersion into $\R^m$ with $L^2$-bounded second fundamental form. Then, up to a subsequence, for any connected component $\sigma$ of $\widetilde{\Sigma}$, the nodal surface of the converging sequence $(\Sigma, \Phi_k^*\xi)$, there exists a M\"obius transformation $\Xi_k$ of $\R^m$ such that 

$$\Xi_k\circ \Phi_k (\Sigma)\subset B(0,R)$$
where $R$ depends only on $m$ and there exists at most finitely many point $\{a_1, \dots, a_L\}$ of $\sigma$ such that if we denote $\Psi_k=\Xi_k \circ \Phi_k\circ\phi_k$, then
$$ \Psi_k \rightharpoonup \Psi \hbox{ weakly in } W^{2,2}_{loc}\left(\sigma \setminus \{a_1, \dots ,a_L, q_1,\dots q_K\} ,\tilde{h}\right),$$
where $\Psi $ is a weak conformal(possibly branched) immersion of $(\sigma,\tilde{h})$ into $\R^m$ and the $q_i$ are the punctures of $(\sigma,\tilde{h})$ and $\phi_k: \widetilde{\Sigma}\rightarrow \Sigma$ such $\phi_k^*(h_k)\rightarrow \tilde{h}$ in $C^\infty_{loc}(\widetilde{\Sigma})$ .\\

Moreover, for any compact $K\subset \sigma \setminus \{a_1, \dots ,a_L, q_1,\dots q_K\} $ there exists $C_K>0$ such that 
$$\sup_{k\in \N} \Vert Log \vert d\Psi_k\vert_{\phi_k^*h_k}\Vert_{L^\infty(K)}\leq C_K ,$$

where $C_K$ depends only on $m$, $K$ and the $L^2$-bound on the  second fundamental form of $\Phi_k$.
\end{thm}

Here we consider the hyperbolic case, since in the sphere case the existence of a non compact conformal group is the additional difficulty already treated in \cite{MondinoRiviere} and in the torus when it degenerate, the injectivity radius uniformly blow down.\\

\noindent{\it Proof of theorem \ref{a2}:}\\

By assumption there exists $\Lambda$, a positive constant , such that 
\beq
\label{hyp}
\sup_{k\in \N }W(\Phi_k) \leq \Lambda .
\eeq

We denote by $u_k$ the conformal factor of this weak immersion with respect to the hyperbolic metric $h_k$ in the conformal class of $\Phi_k^*\xi$. That is to say 
$$ \Phi_k^*\xi=e^{2u_k} h_k ,$$
where $K_{h_k}\equiv -1$.\\

Now let $(\widetilde{\Sigma}, \tilde{h})$ be the nodal surface of the converging sequence $(\Sigma, h_k)$, $\{q_i\}$ the set of punctures, $(\overline{\Sigma}, \overline{h})$ its renormalisation and $\phi_k :\widetilde{\Sigma} \rightarrow \Sigma$ the continuous map given by  to Deligne-Mumford compactification recalled in section 1. Let $\sigma$ be any connected component of $\widetilde{\Sigma}$.\\

Then $\widetilde{\Phi}_k=\Phi_k \circ \phi_k$ is a conformal weak immersion of $(\sigma, \tilde{h}_k)$ where $ \tilde{h}_k =\phi_k^* h_k$. Hence,  we get 
$$ (\Phi_k \circ \phi_k)^*\xi=e^{2u_k} \tilde{h}_k .$$

Let $\delta>0$ and $K_{\delta}=\{x\in\sigma \hbox{ s.t. } d_{\overline{h}}(x, q_i)\geq \delta \hbox{ for all } i\}$, thanks to the local convergence of $\phi_k$, then 
\beq
\label{estu}
\Vert \nabla u_k \Vert_{L^{2,\infty}_{\tilde{h}}(\Sigma)}\hbox{ is uniformly bounded on }K_{\delta}.
\eeq
 Here we use the fact that on the thick part the euclidean metric and the hyperbolic one are equivalent then theorem \ref{a1} can be consider intrinsically on the thick part.\\

Then, in order to find the correct M\"obius transformation, we follow the procedure introduced by the second author in \cite{Riviere10}. For each $x\in K_\delta$ there exists $\rho_{x}^k>0$ such that 
$$\int_{B_{\tilde{h}_k}(x,\rho_x^k)} \vert d\vec{n}_{\widetilde{\Phi}_k} \vert^2_{\tilde{h}_k} dv_{\tilde{h}_k} = \min \left( \frac{8\pi}{3}, \int_{K_\delta} \vert d\vec{n}_{\widetilde{\Phi}_k} \vert^2_{\tilde{h}_k} dv_{\tilde{h}_k}\right) $$
where $B_{\tilde{h}}(x,\rho_x^k)$ is the geodesic ball in $(\sigma, \tilde{h})$ of center $x$ and radius  $\rho_x^k$. Then, using the Besicovitch covering lemma, we can extract a finite covering of  $\ds K_\delta \subset \cup_{ i\in I_k} B_{\tilde{h}}\left(x_i^k,\frac{\rho_{x_i^k}}{2}\right)$, such that each point is covered at most $N$ time, where $N$ is independent of  $k$. Then, thanks to (\ref{hyp}), we can extract a finite covering, $I\subset \cup_k I_k$, which is independent of $k$ where $x_i^k$ converges to $x_i^\infty$ and $\rho_{x_i^k}$ converges to $\rho_{x_i^\infty}$. Then we set $I_0= \{ i\in I \hbox{ s.t. }\rho_{x_i^\infty}=0\}$ and $I_1=I\setminus I_0$.\\

{\bf Claim :For each $i\in I\setminus I_0$ there exist $\overline{v}_k^i\in \R$  such that 
$$\Vert v_k -\overline{v}_k^i\Vert_{L^\infty\left( B_{\tilde{h}}\left(x_i^k,\frac{\rho_{x_i^k}}{2}\right)\right)} \leq C,$$
where $C$ is a constant which depends only on $\Lambda$.}\\

\noindent{\it Proof of the Claim :}\\

Let fix $i\in I\setminus I_0$ and identify, up to uniformly bounded conformal diffeomorphism, $B_{\tilde{h}_k}(x,\rho_x^k)$ with $\D$. Then, thanks to lemma 5.1.4 of \cite{Helein} there exists  a moving frame $(\vec{e}_1^{\, k},\vec{e}_2^{\, k})\in W^{1,2}(\D,S^{m-1})$ such that 
$$\int_{\D}( \vert \nabla \vec{e}_1^{\, k}\vert ^2 +  \vert \nabla \vec{e}_2^{\, k}\vert ^2)\, dz \leq \int_\D \vert d\vec{n}_{\widetilde{\Phi}_k} \vert^2 dz \leq \frac{8\pi}{3},$$
and moreover 
$$  \star \vec{n}_{\widetilde{\Phi}_k} = \vec{e}_1^{\, k} \wedge \vec{e}_2^{\, k} \hbox{ and } \Delta u_k =\left( \nabla^\bot \vec{e}_1^{\, k} , \vec{e}_2^{\, k}\right)  .$$
Let $v_k$ be the solution of 
$$ \left\{\begin{array}{c}\Delta v_k = \left( \nabla^\bot \vec{e}_1^{\, k} , \vec{e}_2^{\, k}\right) \hbox { on } \D \\ v_k = 0 \hbox{ on } \partial \D\end{array}\right. .$$  

Then, thanks to Wente inequality, see section 3 of \cite{Helein}, we get 

\beq
\label{w1}
\Vert v_k \Vert_\infty +\Vert \nabla v_k \Vert_2 \leq \frac{1}{2\pi} \Vert \nabla    \vec{e}_1^k\Vert \Vert \nabla    \vec{e}_2^k\Vert .
\eeq

Finally using the fact $u_k-v_k$ is harmonies with $\Vert \nabla (u_k -v_k)\Vert_{2,\infty}$ uniformly bounded, we proved that there exist $c_k\in \R$ and $C$ a positive constant independent of $k$ such that 

\beq
\label{w2}
 \Vert u_k-v_k- c_k \Vert_{L^\infty \left(B\left(0,\frac{1}{2}\right)\right)} \leq C,
\eeq
 Finally, putting (\ref{w1}) and (\ref{w2}) together concludes the proof of the claim.\hfill$\square$\\

Then using the fact that each point is covered by a universally number of ball, we easily get that there exists $\overline{v}_k\in\R$ such that 

\beq
\label{est2}
\Vert v_k -\overline{v}_k\Vert_{L^\infty\left( K_\delta \setminus \cup_{i\in I_0} B_{\tilde{h}}\left(x_i^\infty,\frac{\delta}{2}\right)\right)} \leq C.
\eeq

We also remark that he constant $\overline{v}_k$ is independent of $\delta$. Let $x_0\in \sigma$ then  we set 

$$\widehat{\Phi}_k =e^{-\overline{v}_k}\left( \tilde{\Phi}_k -\tilde{\Phi}_k(x_0)\right).$$

Then, using  Simon monotonicity formula, see \cite{Simon}, as in \cite{Riviere10} we proved that there exists $y_0\in B(0,1)\subset \R^m$ and $t>0$ such that 
\beq
\label{est3}
\widehat{\Phi}_k \left(K_\delta \setminus \cup_{i\in I_0} B_{\overline{h}}\left(x_i^\infty,\frac{\delta}{2}\right)\right) \cap B(x_0, t)=\emptyset .
\eeq

Finally, we set 
$$\Xi_k = I_{x_0,t} \left( e^{-\overline{v}_k}\left( \ . \ -\tilde{\Phi}_k(x_0)\right)\right)$$
where $I_{x_0,t}$ is the inversion of $\R^m$ centered at $x_0$ and with ratio $t$. Hence $\Xi_k$ is a  M\"obius transformation such that, if we set $\Psi_k= \Xi_k\circ \Phi_k$, we get thanks to (\ref{est2})  that
$$\sup_{k\in \N} \Vert Log \vert d\Psi_k\vert_{\tilde{h}_k}\Vert_{L^\infty(K)}\leq C_K ,$$
and thanks to (\ref{est3}), that  there exist $R>0$ such that,
$$\Psi_k  \left(K_\delta \setminus \cup_{i\in I_0} B_{\overline{h}}\left(x_i,\frac{\delta}{2}\right)\right)\subset B(0,R).$$
Finally, Using a classical argument of functional analysis ,  see for instance \cite{NCRiviere} beginning of section VI.7.1, we easily deduce that $\Psi_k$ converge to $\Psi$ in $W^{2,2}_{loc}\left(\sigma \setminus \{x_1^\infty, \dots ,x_L^\infty, q_1,\dots q_K\} ,\tilde{h}\right)$,
moreover $\Psi$ is a weak immersion away from $ \{x_1^\infty, \dots ,x_L^\infty, q_1,\dots q_K\} $ satisfying
$$W(\Psi) < +\infty$$
Finally  Lemma A.5 of \cite{Riviere11} permits us to extend $\Psi$ as a conformal, possibly branched, immersion of $\sigma$.\hfill$\blacksquare$\\

\section{Weak compactness of Riemannian surfaces with bounded Gaussian curvature entropy}

The last application of our main result is a compactness  result in the spirit of Cheeger and Gromov \cite{CheegerGromov}, Trojanov \cite{Trojanov} and most recently Shioya \cite{Shioya} for Riemannian surfaces. Indeed, we prove a general compactness result for sequence of metrics on a given closed surface assuming only that the area and the total curvature are uniformly bounded  and that the entropy of the Gaussian curvature is also bounded. The first assumptions are the weaker we can assume in order to the  problem makes sense. And the second is made necessarily if one consider a long thin cylinder closed by a two spherical cap, see \cite{Trojanov} and reference therein for more examples of degenerating metrics with bounded curvature and area.

The entropy of the Gaussian curvature of a given metric is defined as follows, let $\Sigma$ be a closed surface and $g$ a Riemannian metric  with Gaussian curvature equal to $K_g$, then we set 
$$E(g)=\int_\Sigma K_g^+ \ln( K_g^+) \, dv_g,$$
where $K_g^+=\max(0,K_g)$ and we set $K_g^+ \ln( K_g^+)=0$ when $K_g^+=0$.
This was introduce by Hamilton in the context of Ricci flow on surfaces. He notably proved that it is monotonically increasing along the Ricci flow on spheres with positive curvature, see \cite{Hamilton} and \cite{Chow}. In order to apply directly our preceding result, we introduce a slightly stronger notion of entropy. Let $\Sigma$ a closed surface with a reference metric $g_0$, then we set 

$$\tilde{E}_0(g)=\int_\Sigma  \ln( e+\vert K_g dv_g\vert_{g_0} ) \vert K_g\vert  \, dv_g.$$

Then considering this notion of entropy, we get the following compactness result.
\begin{thm}
\label{a3} 
For any closed Riemannian surface $(\Sigma, g_0)$ and any  sequence of smooth metric  $g_k$ such that 
$$\int_\Sigma  \vert K_g\vert  \, dv_g + \tilde{E}_0(g_k) =O(1),$$
then for each component $\sigma$ of the thick part of $(\Sigma, h_k)$, then, up to a dilatation of  the metric by a factor $e^{-C_k}$, $h_k$  converges weakly in $(L^{\infty}_{loc}(\sigma))^*$.\\

More precisely, up to a subsequence,  one of the following occurs
\begin{enumerate}[i)]
\item $\mathrm{genus}(\Sigma)=0$, then there exists $C_k$ such that if $e^{-C_k}g_k= e^{2u_k} g_0$, where $g_0$ is the metric of the standard sphere, and $u_k$ is uniformly bounded,
\item  $\mathrm{genus}(\Sigma)=1$, then up to a first dilation, $(\Sigma, g_k)$ is isometric to $\ds \C / \left(\Z \times v_k \Z\right)$ where $v_k$ lies in a fundamental domain of $\Hy/\mathrm{PSL}_2(\Z)$, then there exists $C_k$ such that if  $e^{-C_k}g_k= e^{2u_k} dz^2$ then $u_k$ is  bounded in $L^\infty_{loc}\left( \C / \left(\Z \times v_k \Z\right)\right) $.
\item  $\mathrm{genus}(\Sigma)\geq 1$, then let $\sigma$ a connected component $\sigma$ of the nodal surface of $(\Sigma, h_k)$, then there exists $C_k$ such that if  $e^{-C_k}g_k= e^{2u_k}h_k$ and $u_k$ is bounded in $L^\infty_{loc}(\sigma)$.
\end{enumerate}
\end{thm}

Here, for a sake of simplicity we consider  the standard sphere and a cylinder of fixed radius as thick part.\\

\noindent{\it Proof of theorem \ref{a3}:}\\

 We choose the atlas given by theorem \ref{t1}, and let $g_k= e^{2u_k}h_k$ where $u_k$ is the conformal factor with respect to a normalized metric of constant curvature. Let $K$ be a compact set of $\sigma$ a connected component of the thick part and $U$ an open set of $\sigma$ such that $\overline{U}$ is compact and $K\subset U$. Then on $U$ the conformal factor satisfies
$$\Delta_{h_k} u_k = K_{g_k}e^{2u_k}-K_{h_k},$$ \\
and $h_k$ converges strongly to a smooth metric. On the one hand, since the total curvature is bounded, as in theorem \ref{a1}, we get that $\nabla u_k$ is uniformly bounded in $L^{2,\infty}$.
Then let $v_k \in H^1_0(U)$ such that 
$$\Delta_{h_k} v_k = K_{g_k}e^{2u_k}-K_{h_k} \hbox{ on } U .$$
On the other hand, thanks to the theory of singular integral, see \cite{Stein70} exercise II 6.2.(b),
we get that $\nabla v_k$ is uniformly bounded in $L^2$ and  that $v_k$ is uniformly bounded in $L^{\infty}$. Then, since $u_k-v_k$ is harmonic on $U$ whose gradient is bounded in $L^{2,\infty}$, then thanks to Harnack inequality there exists a constant $C_k$ such that  $u_k -v_k-C_k$ is uniformly bounded on $K$. Then, after checking that the constant is independent of $K$, we get that on each connected component of the thick part there exists a sequence of constant $C_k$ such that $e^{-C_k}g_k= e^{\tilde{u}_k}h_k$ with $\tilde{u}_k$ uniformly bounded in $L^{\infty}_{loc}(\sigma)$, which prove the theorem.\hfill $\square$\\

{\bf Remark :}An interesting problem  is to try to replace $\tilde{E}_0$ by $E$ in the preceding theorem. In order to do  this, on need to analyze the way $K_g dv_g$ concentrates as already done, in the particular case where the Gauss curvature converges uniformly, by \cite{BrezisMerle} see also \cite{LiShafrir}.\\

\end{document}